# Evaluating decision making units under uncertainty using fuzzy multi-objective nonlinear programming

**M. Zerafat Angiz L[1 *], M. K. M. Nawawi[1], R. Khalid[1], A. Mustafa[2], A. Emrouznejad[3], R. John[5], G. Kendall[4, 5]**

[1]*School of Quantitative Sciences, Universiti Utara Malaysia, Sintok, Kedah, Malaysia.*

[2]*School of Mathematical Sciences, Universiti Sains Malaysia, Penang Malaysia.*

[3]*Aston Business School, Aston University, Birmingham, UK.*

[4]*School of Computer Science, University of Nottingham Malaysia Campus, Malaysia*

[5]*ASAP Research Group, School of Computer Science, University of Nottingham, UK*

## Abstract

This paper proposes a new method to evaluate Decision Making Units (DMUs) under uncertainty using fuzzy Data Envelopment Analysis (DEA). In the proposed multi-objective nonlinear programming methodology both the objective functions and the constraints are considered fuzzy. The coefficients of the decision variables in the objective functions and in the constraints, as well as the DMUs under assessment are assumed to be fuzzy numbers with triangular membership functions. A comparison between the current fuzzy DEA models and the proposed method is illustrated by a numerical example.

---

[*] Corresponding author. Tel:+60173509310,

*Email address: madjid@uum.edu.my*

## 1. Introduction

Data Envelopment Analysis (DEA) is a relatively recent approach to the assessment of performance of organizations and their functional units. DEA is able to evaluate the Decision Making Units (DMUs) based on multiple inputs and outputs. Since the first development of DEA [1, 2], there have been many applications of DEA in a variety of different contexts [3, 4].

However in many real applications, input or output variables are not always represented by crisp values. Hence, the traditional DEA models cannot be used for evaluating such DMUs. Several attempts have been made to develop fuzzy DEA models that are powerful tools for comparing the performance of a set of activities or organizations under uncertainty. For instance, Sengupta [5] considered the objective function to be fuzzy when utilizing a standard DEA. He then obtained results using Zimmermann's method [6, 7]. Leon et al. [8] transformed the fuzzy DEA into crisp DEA [9]. Takeda and Satoh [10] used both multicriteria decision analysis and DEA with incomplete data. Lertworasilikul et al. [11, 12] applied the possibilistic approach [13] to treat the constraints of the DEA as fuzzy events. Several other fuzzy models [14] have been proposed to evaluate DMUs with fuzzy data, using the concept of comparison of fuzzy numbers.

The α-cut approach [15, 16] for fuzzy DEA, which is based on a common set of weights in fuzzy DEA, is one of the most frequently used methods. This method first solves a linear program to determine the upper bound of the weights, then a common set of weights are obtained by solving another linear programming problem. We propose a multiobjective programming model [17] that can retain the uncertainty in many aspects including objective functions, coefficients of the decision matrix and the DMUs under assessment.

The paper is organized as follows. A brief description of standard DEA and fuzzy DEA is given in Section 2. A specific multiobjective model will be discussed in section 3. Subsequently, in Section 4 we propose an alternative fuzzy DEA model under uncertainty. This is followed by a numerical illustration in Section 5. In section 6 the methodology is discussed, and mechanism of the proposed approach is presented in section 7. Conclusions are drawn in Section 8.

## 2. DEA and Fuzzy DEA

DEA is a relatively new approach to the assessment of performance of organizations and their functional units. It is a nonparametric technique for measuring the relative efficiency of a set of DMUs with multiple inputs and multiple outputs. Today, DEA is adopted in many disciplines as a powerful tool for assessing efficiency and productivity. Hence many applications of DEA are reported, for example hospital efficiency [18], banking [19, 20], measurement efficiency of health centers [21], manufacturing efficiency [20, 21], productivity of Organisation for Economic Co-operation and Development (OECD) countries [22–24]. Many more applications can be found in the literature [3] which indicates that most of these studies ignored the uncertainty in input and output values. This uncertainty could have an affect on the border defined by the standard DEA; hence the CCR-DEA model may not obtain the true efficiency of DMUs. Theoretically, the standard CCR-DEA model has its production frontier spanned by the linear combination of the observed DMUs.

The production frontier under uncertainty is different. The idea proposed in this research is to allow some flexibility in defining the frontiers with uncertain DMUs, using a fuzzy concept.

### *2.1 Fuzzy number*

**Definition 1.** A triangular fuzzy number $\tilde{x}$ is defined as follows

$$\mu_{\tilde{x}}(\bar{x}) = \begin{cases} \dfrac{\bar{x} - x^l}{x^m - x^l} & \text{for } x^l \le \bar{x} \le x^m \\ \dfrac{x^u - \bar{x}}{x^u - x^m} & \text{for } x^m \le \bar{x} \le x^u \end{cases} \tag{1}$$

$x^m$, $x^l$ and $x^u$ are the mean value, the lower bound and the upper bound of the interval of fuzzy number. The interval of fuzzy number $[x^l, x^u]$ is the region where the value of $\bar{x}$ fluctuates. Symbolically, $\tilde{x}$ is denoted by $(x^m, x^l, x^u)$.

### 2.2 Fuzzy DEA

The technique proposed evaluates the relative efficiency of a set of homogenous DMUs by using a ratio of the weighted sum of outputs to the weighted sum of inputs. It generalizes the usual efficiency measurement from a single-input, single-output ratio to a multiple-input, multiple-output ratio.

Let inputs $x_{ij}$ $(i = 1, 2, ..., m)$ and outputs $y_{rj}$ $(r = 1, 2, ..., s)$ be given for $DMU_j$ $(j = 1, 2, ..., n)$.

The fractional programming statement for the CCR model is formulated as follows:

$$\max \quad \frac{\sum_{r=1}^{s} u_r y_{rp}}{\sum_{i=1}^{m} v_i x_{ip}}$$

$$\text{s.t.}$$

$$\frac{\sum_{r=1}^{s} u_r y_{rj}}{\sum_{i=1}^{m} v_i x_{ij}} \le 1 \qquad \forall j$$

$$u_r, v_i \ge 0 \qquad \forall r, i$$

where $v_i$ and $u_r$ are the weight variables for $i$ th and $r$ th input and output, respectively.

The above model is transformed to the following linear programming problem by some substitutions:

**Model 1: CCR-DEA model**

$$\max \quad \sum_{r=1}^{s} u_r y_{rp}$$

s.t.

$$\sum_{i=1}^{m} v_i x_{ip} = 1$$

$$\sum_{r=1}^{s} u_r y_{rj} - \sum_{i=1}^{m} v_i x_{ij} \leq 0 \qquad \forall j$$

$$u_r, v_i \geq 0 \qquad \forall r, i$$

At the turn of the present century, reducing complex real-world systems into precise mathematical models was the main trend in science and engineering. Unfortunately, real-world situations are frequently not dealing with exact data. Thus precise mathematical models are not enough to tackle all practical problems. In practice there are many problems in which, all (or some) input–output levels are fuzzy numbers. It is difficult to evaluate DMUs in an accurate manner to measure the efficiency. Fuzzy DEA is a powerful tool for evaluating the performance of a set of organizations or activities under an uncertain environment.

Suppose that there are *n* DMUs denoted by *j=1,…,n*, each of which produces a fuzzy nonzero output vector $\tilde{Y}_j = (\tilde{y}_{1j}, \tilde{y}_{2j}, ...\tilde{y}_{sj})^t \geq 0$ using a fuzzy nonzero vector $\tilde{X}_j = (\tilde{x}_{1j}, \tilde{x}_{2j}, ...\tilde{x}_{rj})^t \geq 0$ where the superscript ”t” indicates the transpose of a vector. Consider $\tilde{X}, \tilde{Y}$ are matrices of fuzzy input and output variables of all DMUs. Then, the CCR model with fuzzy coefficients for assessing $DMU_p$ is formulated as follows.

**Model 2: Fuzzy CCR-DEA, multiplier model**

$$\begin{aligned} \max_{u,v} \quad & u^t \tilde{y}_p \\ & v^t \tilde{x}_p = 1 \\ & -v^t \tilde{X} + u^t \tilde{Y} \leq 0 \\ & u^t, v^t \geq 0 \end{aligned}$$

$v^t \in R^{r\times 1}, u^t \in R^{s\times 1}$ column are vectors of inputs and outputs weights, respectively and $\lambda \in R^{n\times 1}$ is column vector of a linear combination of $n$ DMUs.

Saati et al. [15] proposed a fuzzy DEA by considering the α-cut of objective function and the α-cut of constraints; hence the following model is obtained.

**Model 3: Fuzzy CCR-DEA, using α-cut approach**

$$\begin{aligned} \max \quad & \sum_{r=1}^{5} u_r(\alpha y_{rp}^m + (1-\alpha) y_{rp}^l, \alpha y_{rp}^m + (1-\alpha) y_{rp}^u) \\ s.t. \quad & \sum_{i=1}^{m} v_i(\alpha x_{ip}^m + (1-\alpha) x_{ip}^l, \alpha x_{ip}^m + (1-\alpha) x_{ip}^u) = (\alpha + (1-\alpha) l^l, \alpha + (1-\alpha) l^u) \quad \forall_i \\ & \sum_{r=1}^{5} u_r(\alpha y_{rj}^m + (1-\alpha) y_{rj}^l, \alpha y_{rj}^m + (1-\alpha) y_{rj}^u) \\ & -\sum_{i=1}^{m} v_i(\alpha x_{ij}^m + (1-\alpha) x_{ij}^l, \alpha x_{ij}^m + (1-\alpha) x_{ij}^u) \leq 0 \quad \forall j \\ & u_r, v_i \geq 0 \quad \forall r, i. \end{aligned}$$

If we substitute $\tilde{x}_{ij} = (x_{ij}^m, x_{ij}^l, x_{ij}^u)$, $\tilde{y}_{ij} = (y_{ij}^m, y_{ij}^l, y_{ij}^u)$ and $\tilde{1} = (1, 1^l, 1^u)$, Model (3) is written as follows.

**Model 4: Fuzzy CCR-DEA, using α-cut approach, interval programming**

$$\begin{aligned}
\max \quad & \sum_{r=1}^{5} u_r \hat{y}_{rp} \\
s.t. \quad & \sum_{i=1}^{m} v_i \hat{x}_{ip} = L \\
& \sum_{r=1}^{5} u_r \hat{y}_{rj} - \sum_{i=1}^{m} v_i \hat{x}_{ij} \le 0 \\
& \alpha y_{rj}^{m} + (1-\alpha) y_{rj}^{l} \le \hat{y}_{rj} \le \alpha y_{rj}^{m} + (1-\alpha) y_{rj}^{u} \\
& \alpha x_{ij}^{m} + (1-\alpha) x_{ij}^{l} \le \hat{x}_{ij} \le \alpha x_{ij}^{m} + (1-\alpha) x_{ij}^{u} \\
& \alpha + (1-\alpha) l^{l} \le L \le \alpha + (1-\alpha) l^{u} \\
& u_r, v_i \ge 0 \quad \forall r, i.
\end{aligned}$$

As it is shown in Saati et al. [15] we have $\alpha + (1-\alpha)l^{l} \le L \le 1$. One main drawback in Model 4 is that the optimum efficiency level occurs when the outputs of the evaluated DMU and the inputs of other DMUs are set to their upper bounds, while the inputs of the evaluated DMU and the outputs of other DMUs are set to their lower bounds. As a result the evaluated DMU will have the largest possible efficiency value; hence Model 4 does not obtain the true efficiency score.

In the next section we propose an alternative fuzzy DEA to tackle this problem. In the suggested method the evaluated DMU will have the efficiency value between the smallest and the largest possible values.

## 3. Multi--objective programming

Since we must solve a particular multi-objective model, a short discussion related to this kind of problem is presented.

Consider the following multi-objective problem

$$\begin{aligned}
\max \quad & f_1(x), f_2(x), \ldots, f_n(x) \\
\text{s.t.} \quad & \mathrm{x} \in \mathrm{X}
\end{aligned}$$

In the above model, functions $f_1(x), f_2(x),..., f_n(x)$ are objective functions and X is considered as a feasible region. To solve the above mathematical problem, a two stage procedure is proposed.

1. Goal of function $f_i(x)$ $\mathrm{i}=1,2,...,\mathrm{n}$ is obtained by the following mathematical programming:

$$\begin{aligned} & f_i^* = \max\ f_i(x) \\ & \text{s.t.}\quad \mathrm{x} \in \mathrm{X} \end{aligned}$$

2. In this stage scale $\beta$ is introduced to move functions $\frac{f_i(x)}{f_i^*} \le 1$ towards their optimality. For this purpose the following mathematical programming problem should be solved:

$$\begin{aligned} & \max \beta \\ & \text{s.t.}\quad \beta \le \frac{f_i(x)}{f_i^*} \\ & \qquad x \in X \end{aligned}$$

## 3.1. A multi-objective fuzzy DEA model under uncertainty

This section proposes an alternative fuzzy DEA model. The main idea of the suggested method is based on the membership functions of the coefficients. We consider the coefficients as triangular fuzzy numbers $(x^m, x^l, x^u)$. Hence, the membership functions of the coefficients can be defined as follows.

$$\mu_{\tilde{x}_{ij}}(\bar{x}_{ij}) = \begin{cases} \dfrac{\bar{x}_{ij} - x_{ij}^l}{x_{ij}^m - x_{ij}^l} & x_{ij}^l \le \bar{x}_{ij} < x_{ij}^m \\ \dfrac{\bar{x}_{ij} - x_{ij}^u}{x_{ij}^m - x_{ij}^u} & x_{ij}^m \le \bar{x}_{ij} \le x_{ij}^u \end{cases} \qquad \forall i, j \qquad (1)$$

$$\mu_{\tilde{y}_{rj}}(\bar{y}_{rj}) = \begin{cases} \dfrac{\bar{y}_{rj} - y^l_{rj}}{y^m_{rj} - y^l_{rj}} & y^l_{rj} \le \bar{y}_{rj} < y^m_{rj} \\ \dfrac{\bar{y}_{rj} - y^u_{rj}}{y^m_{rj} - y^u_{rj}} & y^m_{rj} \le \bar{y}_{rj} \le y^u_{rj} \end{cases} \quad \forall r, j \qquad (2)$$

Variables $\bar{x}_{ij}$ and $\bar{y}_{rj}$, in formulas (1) and (2), are representative of values in the corresponding intervals of fuzzy numbers.

We suggest the following multi-objective nonlinear program that maximizes both the objective function and the membership functions of technical coefficient simultaneously.

**Model 5: A multi-objective nonlinear programming Fuzzy CCR-DEA**

$$\begin{aligned}
\max \quad & \left\{ \mu_{\tilde{x}_{ij}}(\bar{x}_{ij}), \mu_{\tilde{y}_{rj}}(\bar{y}_{rj}) \right\} \forall j \\
\max \quad & \sum_{r=1}^{s} u_r \bar{y}_{rp} \\
s.t. \quad & \sum_{i=1}^{m} v_i \bar{x}_{ip} = 1 \\
& \sum_{r=1}^{s} u_r \bar{y}_{rj} - \sum_{i=1}^{m} v_i \bar{x}_{ij} \le 0 \quad \forall j (j \ne p) \\
& x^l{}_{ip} \le \bar{x}_{ip} \le x^u{}_{ip} \quad \forall i \\
& y^l{}_{rp} \le \bar{y}_{rp} \le y^u{}_{rp} \quad \forall r \\
& x^l{}_{ij} \le \bar{x}_{ij} \le x^u{}_{ij} \quad \forall i, j \\
& y^l{}_{rj} \le \bar{y}_{rj} \le y^u{}_{rj} \quad \forall r, j \\
& u_r, v_i \ge 0 \quad \forall r, i
\end{aligned}$$

Variables $u_r, v_i$ indicate the coefficients of fuzzy outputs and inputs. Furthermore, variables $\bar{x}_{ij}$ and $\bar{y}_{rj}$ represent the intervals of fuzzy numbers $\tilde{x}_{ij}$ and $\tilde{y}_{rj}$, respectively.

This is a multi-objective nonlinear fuzzy model that we suggest to solve in two stages as explained in the rest of this paper.

Let us ignore the objective functions corresponding to membership functions in Model 5, that is, $\max\left\{\mu_{\tilde{x}_{ij}}(\bar{x}_{ij}), \mu_{\tilde{y}_{rj}}(\bar{y}_{rj})\right\}$ . Then, the optimal solution of the modified model will be as follows:

$$\bar{x}_{ij}^{*} = x_{ij}^{u} \quad j \neq p \quad \bar{x}_{ip}^{*} = x_{ip}^{l}$$
$$\bar{y}_{rj}^{*} = y_{rj}^{l} \quad j \neq p \quad \bar{y}_{rp}^{*} = y_{rp}^{u}$$

This is because each DMU with inputs greater than and outputs less than inputs and outputs $DMU_p$ respectively, will not be better than $DMU_p$. So the optimal value of Model (5) is equals to efficiency of $DMU_p$.

Ignoring the last objective function in Model (5), the optimal solution will be as follows:

$$\bar{x}_{ij}^{*} = x_{ij}^{m} \quad j \neq p \quad \bar{x}_{ip}^{*} = x_{ip}^{m}$$
$$\bar{y}_{rj}^{*} = y_{rj}^{m} \quad j \neq p \quad \bar{y}_{rp}^{*} = y_{rp}^{m}$$

Interaction between two opposed objective functions specify the optimal solution.

**Lemma1:** Let’s consider the optimistic point of view that is the best condition for DMU under evaluation and the worst condition for other DMUs.

a. The optimal solution for $\mu_{\tilde{x}_{ij}}(\bar{x}_{ij}), \mu_{\tilde{y}_{rp}}(\bar{y}_{rp})$ are obtained in the second condition of the membership functions (1) and (2), respectively.
b. The optimal solution for $\mu_{\tilde{x}_{ip}}(\bar{x}_{ip}), \mu_{\tilde{y}'_{rj}}(\bar{y}_{rj})(j \neq p)$ are obtained in the first condition of the membership functions (1) and (2), respectively.

**Proof:** Suppose that objective function in Model (5) be only ( $\max \sum_{r=1}^{s} u_r \bar{y}_{rp}$ ), as mentioned above, due the nature of the model the optimal solution will be:

$$\begin{array}{llll} \min \ \bar{x}_{ip} & \forall i & \max \ \bar{x}_{ij} & \forall i, j(j \neq p) \\ \max \ \bar{y}_{rp} & \forall r & \min \ \bar{y}_{rj} & \forall r, j(j \neq p) \end{array}$$

When considering the effect of the membership function, the values of $\bar{x}_{ij} \quad \forall i, j(j \neq p)$ and $\bar{y}_{rp} \quad \forall r$ will be decreased and the values of $\bar{x}_{ip} \quad \forall i$ and $\bar{y}_{rj} \quad \forall r, j(j \neq p)$ will be increased (membership numbers will be zero for the above mentioned values). So, to obtain the optimal solution of $\mu_{\tilde{x}_{ij}}(\bar{x}_{ij}), \mu_{\tilde{y}_{rp}}(\bar{y}_{rp})$ the second condition of the membership functions (1) and (2) are sufficient, respectively. Similarly to obtain the optimal value for $\mu_{\tilde{x}_{ip}}(\bar{x}_{ip}), \mu_{\tilde{y}'_{rj}}(\bar{y}_{rj})(j \neq p)$ the first condition of the membership functions (1) and (2) are sufficient, respectively, i.e.

$$\mu_{\tilde{x}_{ip}}(\bar{x}_{ip}) = \frac{\bar{x}_{ip} - x^{l}_{ip}}{x^{m}_{ip} - x^{l}_{ip}} \quad \bar{x}_{ip} \in [x^{l}_{ip}, x^{m}_{ip}] \qquad \forall i \tag{3}$$

$$\mu_{\tilde{y}_{rp}}(\bar{y}_{rp}) = \frac{y^{u}{}_{rp} - \bar{y}_{rp}}{y^{u}_{rp} - y^{m}_{rp}} \quad \bar{y}_{rp} \in [y^{m}_{rp}, y^{u}_{rp}] \qquad \forall r \tag{4}$$

$$\mu_{\tilde{x}_{ij}}(\bar{x}_{ij}) = \frac{x^{u}_{ij} - \bar{x}_{ij}}{x^{u}_{ij} - x^{m}_{ij}} \quad \bar{x}_{ij} \in [x^{m}_{ij}, x^{u}_{ij}] \quad \forall i, j(j \neq p) \tag{5}$$

$$\mu_{\tilde{y}_{rj}}(\bar{y}_{rj}) = \frac{\bar{y}_{rj} - y^{l}{}_{rj}}{y^{m}_{rj} - y^{l}_{rj}} \quad \bar{y}_{rj} \in [y^{l}_{rj}, y^{m}_{rj}] \quad \forall r, j(j \neq p) \tag{6}$$

Let $\bar{x}^{*}_{ij}, \bar{y}^{*}_{rj}(j \neq p)$ and $\bar{x}^{*}_{ip}, \bar{y}^{*}_{rp}$ be the optimal solution for $\bar{x}_{ij}, \bar{y}_{rj}(j \neq p)$ and $\bar{x}_{ip}, \bar{y}_{rp}$. It is clear that there exist two values in the intervals $[x^{l}_{ij}, x^{u}_{ij}], [y^{l}_{rj}, y^{u}_{rj}](j \neq p)$ and $[x^{l}_{ip}, x^{u}_{ip}], [y^{l}_{rp}, y^{u}_{rp}]$ with the same membership function, say,

$$\begin{array}{l} \bar{x}^{*}_{ij1} \in [x^{l}_{ij}, x^{m}_{ij}], \bar{y}^{*}_{rj1} \in [y^{l}_{rj}, y^{m}_{rj}] \\ \bar{x}^{*}_{ij2} \in [x^{m}_{ij}, x^{u}_{ij}], \bar{y}^{*}_{rj2} \in [y^{m}_{rj}, y^{u}_{rj}] \end{array} \tag{7}$$

$$\bar{x}^*_{ip1} \in [x^l_{ip}, x^m_{ip}], \bar{y}^*_{rp1} \in [y^l_{rp}, y^m_{rp}]$$

$$\bar{x}^*_{ip2} \in [x^m_{ip}, x^u_{ip}], \bar{y}^*_{rp2} \in [y^m_{rp}, y^u_{rp}].$$

In this view, the $\bar{x}_{ij}s$ are similar to the input values and the $\bar{y}_{rj}s$ are similar to the output values in the DEA models, so by considering constant values for $\bar{x}_{ij}s$ and $\bar{y}_{rj}s$, Model (5) will be converted to Model (4).

Assume that inputs and outputs of $DMU_1$ and $DMU_2$ are $(x'^*_{ip1}, x'^*_{ij1}, y'^*_{rp2}, y'^*_{rj2})(j \neq p)$ and $(x'^*_{ip2} x'^*_{ij2}, y'^*_{rp1}, y'^*_{rj1})(j \neq p)$, respectively. Obviously $DMU_1$ is more efficient than $DMU_2$. This means only the second condition of the membership functions (1) and (2) are sufficient to obtain the optimal solution for $\mu_{\tilde{x}'_{ij}}(x'_{ij}), \mu_{\tilde{y}'_{ip}}(y'_{ip})$. Similarly the first condition of membership function (1) and (2) are sufficient to obtain the optimum value for $\mu_{\tilde{x}'_{ip}}(x'_{ip}), \mu_{\tilde{y}'_{ij}}(y'_{ij})(j \neq p)$.

Hence, to solve Model (5), the methodology presented in section 3 is applied, and multi-objective programming problem (5) is converted to the following nonlinear programming problem:

**Model 6: A new Fuzzy CCR-DEA, non-linear programming**

$$\max \quad Z = h$$

$$s.t.$$

$$\sum_{i=1}^{m} v_i \bar{x}_{ip} = 1$$

$$h \le \left(\sum_{r=1}^{s} u_r \bar{y}_{rp}\right) / z_p^{*}$$

$$\sum_{r=1}^{s} u_r \bar{y}_{rj} - \sum_{i=1}^{m} v_i \bar{x}_{ij} \le 0 \quad \forall j (j \ne p)$$

$$h \le \frac{x^{u}{}_{ij} - \bar{x}_{ij}}{x^{u}_{ij} - x^{m}_{ij}} \qquad \forall i, j (j \ne p)$$

$$h \le \frac{\bar{y}_{rj} - y^{l}_{rj}}{y^{m}_{rj} - y^{l}_{rj}} \qquad \forall r, j (j \ne p)$$

$$h \le \frac{\bar{x}_{ip} - x^{l}_{ip}}{x^{m}_{ip} - x^{l}_{ip}} \qquad \forall i$$

$$h \le \frac{y^{u}{}_{rp} - y_{rp}}{y^{u}_{rp} - y^{m}_{rp}} \qquad \forall r$$

$$x^{m}_{ij} \le \bar{x}_{ij} \le x^{u}_{ij} \qquad \forall i, j (j \ne p) \qquad 6.1$$

$$y^{l}_{rj} \le \bar{y}_{rj} \le y^{m}_{rj} \qquad \forall r, j (j \ne p) \qquad 6.2$$

$$x^{l}_{ip} \le \bar{x}_{ip} \le x^{m}_{ip} \qquad \forall i \qquad 6.3$$

$$y^{m}_{rp} \le \bar{y}_{rp} \le y^{u}_{rp} \qquad \forall r \qquad 6.4$$

$$u_r, v_i \ge 0 \quad \forall r, i$$

In Model (6), $z_p^*$ is obtained with the best situation of the DMUs as follows:

**Model 7: A new Fuzzy CCR-DEA, estimation of $Z^*_p$**

$$z_p = \max \quad \sum_{r=1}^{5} u_r y^{u}_{rp}$$

$$s.t \qquad \sum_{i=1}^{m} v_i x^{l}_{ip} = 1$$

$$\sum_{r=1}^{s} u_r y^{l}_{rj} - \sum_{i=1}^{m} v_i x^{l}_{ip} \le 0 \quad \forall j (j \ne p)$$

$$u_r, v_i \ge 0 \quad \forall r, i$$

Obviously, fluctuating between 0 and 1, the objective functions corresponding to membership functions do not need to follow the first stage of section 3.

The variable *h* in Model (6) is used to convert the multi-objective problem Model (5) to a nonlinear programming problem. This variable is within the interval $[0,1]$. Adding the concept of α-cut to Model (6), it is sufficient to replace the following constraints instead of 6-1, 6-2, 6-3 and 6-4.

$$
\begin{array}{ll}
x'^{m}_{ij} \le x'_{ij} \le \alpha x'^{m}_{ij} + (1-\alpha) x'^{u}_{ij} & \forall i, j (j \ne p) \\
\alpha y'^{m}_{rj} + (1-\alpha) y'^{l}_{rj} \le y'_{rj} \le y'^{m}_{rj} & \forall r, j (j \ne p) \\
\alpha x'^{m}_{ip} + (1-\alpha) x'^{l}_{ip} \le x'_{ip} \le x'^{m}_{ip} & \forall i \\
y'^{m}_{rp} \le y'_{rp} \le \alpha y'^{m}_{rp} + (1-\alpha) y'^{u}_{rp} & \forall r
\end{array}
$$

This is different from the standard α-cut used in the fuzzy DEA Model (4), because in each α-level the model still retains uncertainty information interior of the interval that was generated by α. Next section compares our results with the current fuzzy DEA model.

## 4. An illustration with a numerical example

In this section, a numerical example is presented to illustrate the difference between the results obtained using the proposed approach and the current fuzzy DEA models. Consider the data in Table 1 that is extracted from Guo and Tanaka [14] and used by Lertworasirikul et al. [11] and Saati et al. [15]. There are 5 DMUs with two symmetrical triangular fuzzy inputs and 2 symmetrical triangular fuzzy outputs.

**Table 1: Data for numerical example**

| Variable | DMU D1 | D2 | D3 | D4 | D5 |
|---|---|---|---|---|---|
| I1 | (4.0, 3.5, 4.5) | (2.9, 2.9, 2.9) | (4.9, 4.4, 5.4) | (4.1, 3.4, 4.8) | (6.5, 5.9, 7.1) |
| I2 | (2.1, 1.9, 2.3) | (1.5, 1.4, 1.6) | (2.6, 2.2, 3.0) | (2.3, 2.2, 2.4) | (4.1, 3.6, 4.6) |
| O1 | (2.6, 2.4, 2.8 | (2.2, 2.2, 2.2) | (3.2, 2.7, 3.7) | (2.9, 2.5, 3..3) | (5.1, 4.4, 5.8) |
| O2 | (4.1, 3.8, 4.4) | (3.5, 3.3, 3.7) | (5.1, 4.3, 5.9) | (5.7, 5.5, 5.9) | (7.4, 6.5, 8.3) |

Using fuzzy CCR Model (4), the efficiency scores are summarized in the Table 2.

**Table 2: The efficiencies using Model (4)**

| A | DMU D1 | D2 | D3 | D4 | D5 |
|---|---|---|---|---|---|
| 0 | 1.11 | 1.51 | 1.28 | 1.52 | 1.30 |
| .5 | .995 | 1.32 | 1.03 | 1.32 | 1.16 |
| .75 | .906 | 1.24 | 0.93 | 1.23 | 1.12 |
| 1 | .85 | 1 | .86 | 1 | 1 |

Considering the above Lemma1, obviously, the optimal solution given in Table 2 is equivalent to the optimal solution related to the optimistic part of Kao and Liu [27] approach in its supper efficiency form. As it is known, the methods based on the **α**-cut approach just extent number of membership values considered in the evaluation; therefore the major part of the fuzzy concept is ignored. Differences between the proposed method and the **α**-cut based approach can be compared with differences between integration and numerical methods for integrals. The numerical methods don’t cover the whole area under curve in integration.

Results from the possibility approach of Lertworasirikul [11] are shown in Table 3. As can be seen, the efficiency values in the above two models are very similar.

**Table 3: The efficiencies using Lertworasirikul** [11] **model**

| α | DMU | | | | |
|---|---|---|---|---|---|
| | D1 | D2 | D3 | D4 | D5 |
| 0 | 1.11 | 1.24 | 1.28 | 1.52 | 1.30 |
| .5 | 0.96 | 1.11 | 1.03 | 1.26 | 1.16 |
| .75 | .91 | 1.06 | 0.93 | 1.13 | 1.10 |
| 1 | .85 | 1 | .86 | 1 | 1 |

Using the proposed Model (6), the results are shown in Table 4.

**Table 4: The efficiencies using the proposed model in this paper**

| α | DMU | | | | |
|---|---|---|---|---|---|
| | D1 | D2 | D3 | D4 | D5 |
| 0 | 0.899 | 1.220 | 0.930 | 1.220 | 1.076 |
| 0.5 | 0.86 | 1.180 | 0.871 | 1.169 | 1.041 |
| 0.75 | 0.85 | 1.110 | 0.866 | 1.160 | 1.037 |
| 1 | 0.85 | 1.000 | 0.860 | 1.000 | 1.000 |

Due to the nature of the fuzzy CCR Model (4) the maximum efficiency occurs when the outputs of the evaluated DMU and the inputs of other DMUs are set to their upper bounds. It is obvious that the results in Table 2 are always greater than the results that we obtained in Table 4 since Model 4 always captures the efficiency under pessimistic circumstances. The results obtained using the proposed model in this paper have the efficiency values between the smallest and the largest possible values, hence they are more close to the true efficiency.

## 6. Empirical study

To illustrate the fuzzy DEA approach, we consider data given in [28] which has presented for an aircraft selection. Five types of aircraft (B757-200, A-321, B767-200, MD-82, and A310-300) are to be evaluated. Four inputs and two outputs are introduced in Table 5 as follows:

**Table 5: Inputs and outputs for aircrafts evaluation**

| Data | Description |
|---|---|
| Input1(I1) | Maintenance requirements (Subjective assessment) |
| Input2(I2) | Pilot adaptability (Subjective assessment) |
| Input3(I3) | Maximum range (Kilometer) |
| Input4(I4) | Purchasing price (US millions) |
| Output1(O1) | Passenger preference(Subjective assessment) |
| Output2(O2) | Operational productivity (Seat-kilometer per hour) |

The first input is the aircraft maintenance capability (I1) which is concerned with the availability and the level of standardization of spare parts and post-sale services. The second input, pilot adaptability (I2) is related to the skills of available pilots and the specific features of the aircraft. The third input maximum range (I3) of an aircraft is determined by the maximum kilometers that the aircraft can travel at the maximum payload and the fourth input, purchasing price (I4) is the price to be paid for a new aircraft which correlates with reliability of the aircraft.

On the other hand for the outputs, passengers' preference (O1) reflects the social responsibility of the airline in order to establish a positive image in public and of the requirements imposed by various environment protection laws and regulations whilst operational productivity (O2) is determined by the number of seats available, the load rate, the travel frequency, and the aircraft travel speed.

In this research, the eight decision makers stated their opinion about 3 subjective inputs and outputs. They used a set of five linguistic terms {very low, low, medium, high, very high} which are associated with the corresponding numbers 1, 2, 3, 4 and 5, respectively, as in a 5-point Likert scale.

Table 6 shows the inputs and outputs of the five aircrafts. For example, B757-200 type of aircraft has two subjective inputs (I1 and I2) and one subjective output (O1), with triangular fuzzy numbers. For other two inputs and one output, the values are crisps.

**Table 6: Data for numerical example**

| | DMU | | | | |
|---|---|---|---|---|---|
| Variable | B757-200 | A-321 | B767-200 | MD-82 | A310-300 |
| I1 | (2.0, 3.064, 4) | (4, 4.229,5) | (3, 3.224, 4) | (1, 1.929, 3) | (3,3.464, 4) |
| I2 | (2, 2.852, 3) | (2,2.000,2) | (2, 2.852, 3) | (4, 4.113, 5) | (2,2.000,2) |
| I3 | 5522 | 4350 | 5856 | 4032 | 7968 |
| I4 | 56 | 54 | 69 | 33 | 80 |
| O1 | (4, 4.000, 4) | (2, 2.852, 3) | (4, 4.000, 4) | (3, 3.591, 4) | (3, 3.342, 4) |
| O2 | 116279 | 109063 | 129465 | 87662 | 130664 |

Using Model (6), the values of h*, the efficiency scores and rank of each aircraft are given in Table 7. The MD-82 aircraft type gives the highest efficiency score of 1.8520 and is ranked first, whilst B767-200 gives lowest score of 1.0949 and is ranked last.

**Table 7: The rank of five types of aircrafts**

| DMU | $h^*$ | Eff. scores | Rank |
|---|---|---|---|
| B757-200 | 0.6348 | 1.2696 | 2 |
| A-321 | 0.9798 | 1.1720 | 3 |
| B767-200 | 1.0000 | 1.0949 | 5 |
| MD-82 | 0.9260 | 1.8520 | 1 |
| A310-300 | 1.0000 | 1.1237 | 4 |

## 5. Discussion

According to theorem 2, if the objective functions corresponding to membership functions in Model (5) are ignored, the optimal solution for inputs and outputs will be arisen in endpoints of interval of fuzzy numbers. Furthermore, if the last objective function ($\max \sum_{r=1}^{s} u_r \overline{y}_{rp}$ ) in Model (5) is eliminated, Lemma1 adopted the optimal solution will be in the main value of fuzzy number.  Figure 1 illustrates the above mentioned concept for evaluating $DMU_P$.  The interior arrows represent the optimal solution when the last objective function ($\max \sum_{r=1}^{s} u_r \overline{y}_{rp}$ ) is absent in Model (5) and the arrows located under fuzzy numbers construct the optimal solution Model (5) when only the objective function ($\max \sum_{r=1}^{s} u_r \overline{y}_{rp}$) is present.

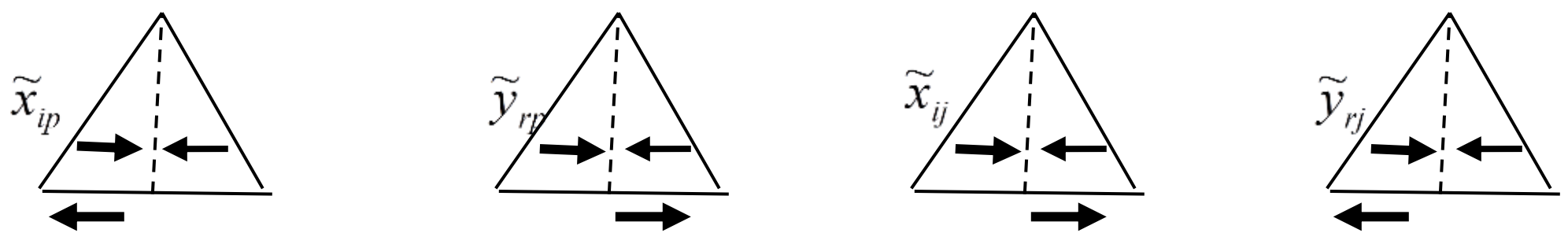

**Figure 1: Concepts of evaluating DMUs**

Interaction between the objective functions corresponding to objective functions and the last objective function ($\max \sum_{r=1}^{s} u_r \bar{y}_{rp}$) in Model (5), cause the fuzzy optimal solution.

## 6. Conclusion

In evaluating DMUs under uncertainty several fuzzy DEA models have been proposed in the literature. The **α**-cut approach is one of the most frequently used models. However, due to the nature of the **α**-cut approach the uncertainty in inputs and outputs is effectively ignored. This paper proposed a multi-objective fuzzy DEA model to retain fuzziness of the model by maximizing the membership function of inputs and outputs. In the proposed method, both the objective functions and the constraints are considered fuzzy. A numerical example is used to show the difference between the proposed and the current fuzzy DEA models. For further studies, it is suggested that an exploration be done on: a) reducing the size of the converted (crisp equivalent) problem, b) possible linearization of the nonlinear model.